\documentclass[12pt]{article}
\usepackage{amsfonts}
\usepackage{mathrsfs}
\usepackage{amsmath,amsfonts,amssymb,rotating}

\newtheorem{theo}{Theorem}[section]

\newtheorem{conj}[theo]{Conjecture}

\newdimen\Squaresize \Squaresize=11pt
\newdimen\Thickness \Thickness=0.7pt
\def\Square#1{\hbox{\vrule width \Thickness
   \vbox to \Squaresize{\hrule height \Thickness\vss
    \hbox to \Squaresize{\hss#1\hss}
   \vss\hrule height\Thickness}
\unskip\vrule width \Thickness} \kern-\Thickness}

\def\Vsquare#1{\vbox{\Square{$#1$}}\kern-\Thickness}

\def\moins{\raise 1pt\hbox{{$\scriptstyle -$}}}

\begin{document}

\begin{center}
{\large \bf Log-concavity and $q$-Log-convexity Conjectures on the

Longest Increasing Subsequences of Permutations}
\end{center}

\begin{center}
William Y.C. Chen\\[6pt]
Center for Combinatorics, LPMC-TJKLC\\
Nankai University, Tianjin 300071, P. R. China

  ${\tt
chen@nankai.edu.cn}$
\end{center}

\vspace{0.2cm} \noindent{\bf Abstract.} Let $P_{n,k}$ be
the number of permutations $\pi$ on $[n]=\{1, 2, \ldots,
n\}$ such that the length of the longest increasing
subsequences of $\pi$ equals $k$, and let $M_{2n, k}$ be
the number of matchings  on $[2n]$ with crossing number
$k$. Define $P_n(x)= \sum_k P_{n,k}x^k$ and
$M_{2n}(x)=\sum_{k} M_{2n,k}x^k$. We propose some
conjectures on the log-concavity and $q$-log-convexity of
the polynomials $P_n(x)$ and $M_{2n}(x)$. We also introduce
the notions of $\infty$-$q$-log-convexity and
$\infty$-$q$-log-concavity, and the notion of higher order
log-concavity with respect to $\infty$-$q$-log-convex or
$\infty$-$q$-log-concavity. A conjecture on the
$\infty$-$q$-log-convexity of the Boros-Moll polynomials is
presented. It seems that $M_{2n}(x)$ are log-concave of any
order with respect to $\infty$-$q$-log-convexity.



\section{The Conjectures}

Let $P_n(x)$ and $M_{2n}(x)$ be defined as in the abstract. We
propose the following conjectures.

\begin{conj}\label{con1}
 $P_n(x)$ is log-concave for $n\geq 1$.
\end{conj}

\begin{conj}\label{con2}
 $P_n(x)$ is $\infty$-log-concave for $n\geq 1$.
\end{conj}

\begin{conj}\label{con3}
The polynomial sequence $\{P_n(x)\}$ is strongly $q$-log-convex.
\end{conj}

\begin{conj}\label{con4}
The polynomial sequence $\{P_n(x)\}$ is $\infty$-$q$-log-convex.
\end{conj}

\begin{conj}\label{con5}
 $M_{2n}(x)$ is log-concave for $n\geq 1$.
\end{conj}

\begin{conj}\label{con6}
 $M_{2n}(x)$ is $\infty$-log-concave for $n\geq 1$.
\end{conj}

\begin{conj}\label{con7}
The polynomial sequence $\{M_{2n}(x)\}$ is strongly $q$-log-convex.
\end{conj}

\begin{conj}\label{conj8}
The polynomial sequence $\{M_{2n}(x)\}$ is $\infty$-$q$-log-convex.
Furthermore, the polynomials $M_{2n}(x)$ are log-concave of any
order with respect to $\infty$-$q$-log-convexity.
\end{conj}

The following conjecture is concerned with the Boros-Moll
polynomials¡¡ \cite{BM, BM2}. The log-concavity is established by
Kauser and Paule \cite{kau}.

\begin{conj} The sequence of the  Boros-Moll polynomials is
$\infty$-$q$-log-convex, and they   are log-concave of any order
with respect to $\infty$-$q$-log-convexity.
\end{conj}

\section{The Background}

The longest increasing subsequences of permutations have been
extensively studied; see, for example, \cite{baik1,Berg1, chen,
gessel, regev}, in particular, the survey of Stanley \cite{stan4}.
Baik, Deift and Johansson \cite{baik1} have shown that the limiting
distribution of the coefficients of $P_n(x)$ is the Tracy-Widom
distribution. The numbers $P_{n,k}$ can be computed by Gessel's
theorem \cite{gessel}. Let $\mathfrak{S}_n$ be the symmetric group
on $[n]$, and let ${\rm is}(\pi)$ be the length of the longest
increasing subsequences of $\pi$. Define
\begin{eqnarray}
u_k(n) & = & \sharp \{w \in \mathfrak{S}_n : {\rm is}(w)\leq
k\}, \\[8pt]
U_k(x) & = & \sum \limits_{n\geq 0} u_k(n)\frac{x^{2n}}{n!^2},\quad
k\geq 1, \\[8pt]
I_i(2x) & = & \sum \limits_{n\geq 0} \frac{x^{2n+i}}{n!(n+i)!},\quad
i\in \mathbb{Z}.
\end{eqnarray}

\begin{theo}
\begin{equation}
U_k(x)=\det(I_{i-j}(2x))_{i,j=1}^k.
\end{equation}
\end{theo}

Since $P_{n,k}=u_k(n)-u_{k-1}(n)$ for $n\geq 1$, we can use Gessel's
theorem to compute $P_{n,k}$ for small $n$. Here we list $P_n(x)$
for $1\leq n\leq 18$: \allowdisplaybreaks
\begin{align*}
&P_1(x)= x,\\[5pt]
&P_2(x)= x+x^2,\\[5pt]
&P_3(x)= x+4x^2+x^3,\\[5pt]
&P_4(x)= x+13x^2+9x^3+x^4,\\[5pt]
&P_5(x)= x+41x^2+61x^3+16x^4+x^5,\\[5pt]
&P_6(x)= x+131x^2+381x^3+181x^4+25x^5+x^6,\\[5pt]
&P_7(x)= x+428x^2+2332x^3+1821x^4+421x^5+36x^6+x^7,\\[5pt]
&P_8(x)= x+1429x^2+14337x^3+17557x^4+6105x^5+841x^6+49x^7+x^8\\[5pt]
&P_9(x)=x+4861x^2+89497x^3+167449x^4+83029x^5+16465x^6+1513x^7\\[5pt]
& \quad \quad \quad \quad+64x^8+x^9,\\[5pt]
&P_{10}(x)=x+16795x^2+569794x^3+1604098x^4+1100902x^5+296326x^6\\[5pt]
&  \quad \quad \quad \quad+38281x^7+2521x^8+81x^9+x^{10}.\\[5pt]
&P_{11}(x)=x+58785x^2+3704504x^3+15555398x^4+14516426x^5+5122877x^6\\[5pt]
&  \quad \quad \quad \quad+874886x^7+79861x^8+3961x^9+100x^{10}+x^{11}\\[5pt]
&P_{12}(x)=x+208011x^2+24584693x^3+153315999x^4+192422979x^5+87116283x^6\\[5pt]
&  \quad \quad \quad
\quad+18943343x^7+2250887x^8+153341x^9+5941x^{10}+121x^{11}+x^{12}\\[5pt]
&P_{13}(x)=x+742899x^2+166335677x^3+1538907306x^4+2579725656x^5\\[5pt]
&  \quad \quad \quad \quad+1477313976x^6+399080475x^7+59367101x^8+5213287x^9+275705x^{10}\\[5pt]
&  \quad \quad \quad \quad+8581x^{11}+144x^{12}+x^{13}\\[5pt]
&P_{14}(x)=x+2674439x^2+1145533650x^3+15743413076x^4+35098717902x^5\\[5pt]
&  \quad \quad \quad \quad+25191909848x^6+8312317976x^7+1508071384x^8+164060352x^9\\[5pt]
&  \quad \quad \quad \quad+11110464x^{10}+469925x^{11}+12013x^{12}+169x^{13}+x^{14}\\[5pt]
&P_{15}(x)=x+9694844x^2+8017098273x^3+164161815768x^4+485534447114x^5\\[5pt]
&  \quad \quad \quad \quad+434119587475x^6+172912977525x^7+37558353900x^8\\[5pt]
&  \quad \quad \quad \quad+4927007100x^9+410474625x^{10}+22128576x^{11}+766221x^{12}\\[5pt]
&  \quad \quad \quad \quad+16381x^{13}+196x^{14}+x^{15}\\[5pt]
&P_{16}(x)=x+35357669x^2+56928364553x^3+1744049683213x^4\\[5pt]
&  \quad \quad \quad \quad+6835409506841x^5+7583461369373x^6+3615907795025x^7\\[5pt]
&  \quad \quad \quad \quad+927716186325x^8+143938455225x^9+14353045401x^{10}+947236425x^{11}\\[5pt]
&  \quad \quad \quad \quad+41662441x^{12}+1203441x^{13}+21841x^{14}+225x^{15}+x^{16}\\[5pt]
&P_{17}(x)=x+129644789x^2+409558170361x^3+18865209953045x^4\\[5pt]
&  \quad \quad \quad
\quad+97966603326993x^5+134533482045389x^6+76340522760097x^7\\[5pt]
&  \quad \quad \quad \quad+22904111472825x^8+4142947526101x^9+484748595081x^{10}\\[5pt]
&  \quad \quad \quad
\quad+38094121561x^{11}+2043822961x^{12}+74797417x^{13}+1830561x^{14}\\[5pt]
&  \quad \quad \quad \quad+28561x^{15}+256x^{16}+x^{17}.\\[5pt]
&P_{18}(x)=x+477638699x^2+2981386305018x^3+207591285198178x^4\\[5pt]
&  \quad \quad \quad
\quad+1429401763567226x^5+2426299018270338x^6+1631788075873114x^7\\[5pt]
&  \quad \quad \quad
\quad+568209449266202x^8+118504614869214x^9+16029615164446x^{10}\\[5pt]
&  \quad \quad \quad
\quad+1470147102730x^{11}+93574631242x^{12}+4166173834x^{13}\\[5pt]
&  \quad \quad \quad
\quad+128922442x^{14}+2708305x^{15}+36721x^{16}+289x^{17}+x^{18}.
\end{align*}

One can check that $P_n(x)$ are log-concave for $1\leq n\leq 18$. We
now recall the notion of $k$-log-concavity; see, \cite{kau}. Define
the operator $\mathcal {L}$ which maps a sequence $\{a_i\}$ of
nonnegative numbers to a sequence $\{b_i\}$ given by
$$b_i:= a_i^2-a_{i-1}a_{i+1}.$$ Then the log-concavity of the
sequence $\{a_i\}$ is defined by the positivity of
$\mathcal{L}\{a_i\}$, namely, $b_i$ is nonnegative  for all $i$. If
the sequence $\mathcal{L}\{a_i\}$ is not only positive but also
log-concave, then we say that $\{a_i\}$ is $2$-log-concave. In
general, we say that $\{a_i\}$ is $k$-log-concave if
$\mathcal{L}^k\{a_i\}$ is nonnegative, and that $\{a_i\}$ is
$\infty$-log-concave if $\mathcal{L}^k\{a_i\}$ is nonnegative for
every $k\geq 1$.

In fact, when $n\leq 18$, we can find the sequence
$\{P_{n,k}\}_{1\leq k\leq n}$ is $4$-log-concave. This evidence
leads us to surmise  that the sequence $\{P_{n,k}\}_{1\leq k\leq n}$
is $\infty$-log-concave.

The $q$-log-concavity of polynomials  has  been  studied by many
authors including Butler \cite{butler1990}, Krattenthaler
\cite{kratte1989}, Leroux \cite{leroux1990}, and Sagan
\cite{sagan1992,sagan1992t}. Notice that here we have use $x$
instead of $q$ for the polynomials $P_n(x)$ and $M_{2n}(x)$.
Following the notation of Sagan \cite{sagan1992t}, given two
polynomials $f(q)$ and $g(q)$ in $q$, we write
$$f(q)\geq_q g(q)$$ if the difference $f(q)-g(q)$ has nonnegative
coefficients as a polynomial of $q$. A sequence of polynomials
$\{f_k(q)\}_{k\geq 0}$ over the field of real numbers is called
$q$-log-concave if
\begin{equation}
f_{m}(q)^2\geq_q f_{m+1}(q)f_{m-1}(q),\quad \mbox{ for all $m\geq
1$}.
\end{equation}
 Liu and
Wang \cite{lw} introduced the notion of  $q$-log-convexity. A
polynomial sequence $\{f_n(q)\}_{n\geq 0}$ is called  $q$-log-convex
if
\begin{equation}
f_{m+1}(q)f_{m-1}(q)\geq_q f_{m}(q)^2,\quad \mbox{ for all $m\geq
1$}.
\end{equation}
 A stronger property, called strong
$q$-log-convexity, is introduced by Chen, Wang and Yang
\cite{chen2}. A polynomial sequence  $\{f_n(q)\}_{n\geq 1}$ is
called strongly $q$-log-convex if
\begin{equation}
f_{m-1}(q)f_{n+1}(q)\geq_q f_{m}(q)f_{n}(q),\quad \mbox{ for all
$1\leq m\leq n$}.
\end{equation}
 When $1\leq n\leq 17$, we find
$P_{m-1}(x)P_{n+1}(x)\geq_x P_{n}(x)P_{m}(x)$.

Motivated by the notion of  $\infty$-log-concavity, we define the
operator $\mathcal {H}$ which maps a polynomial sequence
$\{A_i(q)\}_{i\geq 0}$ to a polynomial sequence $\{B_i(q)\}_{i\geq
0}$ given by
$$B_i(q):=A_{i-1}(q)A_{i+1}(q)-A_i(q)^2.$$ Then the $q$-log-convexity of the
polynomial sequence $\{A_i(q)\}$ is defined by the $q$-positivity of
$\mathcal{H}\{A_i(q)\}$, namely, the coefficients of $B_i$ are
nonnegative for all $i$. If the polynomial sequence $\{B_i(q)\}$ is
$q$-log-convex, then we say that $\{A_i(q)\}$ is $2$-$q$-log-convex.
In general, we say that $\{A_i(q)\}$ is $k$-$q$-log-convex if the
coefficients of $\mathcal{H}^k\{A_i(q)\}$ are nonnegative, and that
$\{A_i(q)\}$ is $\infty$-$q$-log-convex if $\mathcal{H}^k\{A_i(q)\}$
is nonnegative for every $k\geq 1$.

 When
$1\leq n\leq 16$, we find the polynomial sequence $\{P_n(x)\}$
log-concave of order 3 with respect to  $3$-$q$-log-convex. This
leads us to surmise the polynomial sequence $\{P_n(x)\}$ is
$\infty$-$q$-log-convex.

We now give a brief review on how to compute the polynomials
$M_{2n}(x)$. The crossing number of a matching on $[2n]$ is the
maximum number $k$ such that there are $k$ mutually intersecting
edges in the standard representation of the matching; see
\cite{chen}. Let $v_{k}(n)$ denote the number of matchings on $[2n]$
whose crossing number is not greater than $k$. For example, the
crossing number of a noncrossing matching equals one, since by
``noncrossing'' we really mean $2$-noncrossing. Note that here we
have used a slightly different notation from that in \cite{stan4}.
 Let
\begin{equation}
V_{k}(x)=\sum \limits_{n\geq 0}v_{k}(n)\frac{x^n}{n!}.
\end{equation}
Grabiner and Magyar \cite{gra} derived the following matching
analogue of Gessel's Theorem. The same formula has also been
obtained by Goulden \cite{gou}.

\begin{theo}
\begin{equation}
V_{k}(x)=\det (I_{i-j}(2x)-I_{i+j}(2x))_{i,j=1}^k.
\end{equation}
\end{theo}

Applying the above theorem, we can compute $v_k(n)$ when $n$ is
small. Since $M_{2n,k}=v_{k}(n)-v_{k-1}(n)$, we obtain
\allowdisplaybreaks
\begin{align*}
&M_2(x)=x\\[5pt]
&M_4(x)=2x+x^2\\[5pt]
&M_6(x)=5x+9x^2+x^3\\[5pt]
&M_8(x)=14x+70x^2+20x^3+x^4\\[5pt]
&M_{10}(x)=42x+552x^2+315x^3+35x^4+x^5\\[5pt]
&M_{12}(x)=132x+4587x^2+4730x^3+891x^4+54x^5+x^6\\[5pt]
&M_{14}(x)=429x+40469x^2+71500x^3+20657x^4+2002x^5+77x^6+x^7\\[5pt]
&M_{16}(x)=1430x+377806x^2+1110174x^3+468650x^4+64960x^5\\[5pt]
&\quad \quad \quad \quad+3900x^6+104x^7+x^8\\[5pt]
&M_{18}(x)=4862x+3707054x^2+17850170x^3+10717004x^4+2005830x^5\\[5pt]
&\quad \quad \quad \quad+167484 x^6+6885 x^7+135 x^8+x^9\\[5pt]
&M_{20}(x)=16796x+37958960x^2+298110266x^3+250367036x^4+61205916x^5\\[5pt]
&\quad \quad \quad
\quad+6681255x^6+376770x^7+11305x^8+170x^9+x^{10}\\[5pt]
&M_{22}(x)=58786x+403068470x^2+5174115036x^3+6012729626x^4+1881276355x^5\\[5pt]
&\quad \quad \quad
\quad+258507711x^6+18770290x^7+766535x^8+17556x^9+209x^{10}+x^{11}\\[5pt]
&M_{24}(x)=208012x+4414995268x^2+93255969556x^3+148847198843x^4\\[5pt]
&\quad \quad \quad
\quad+58846367560x^5+9929079622x^6+892328976x^7+46525941x^8\\[5pt]
&\quad \quad \quad \quad+1443112x^9+26082x^{10}+252x^{11}+x^{12}\\[5pt]
&M_{26}(x)=742900x+49670294000x^2+1742677176125x^3+3801821241675x^4\\[5pt]
&\quad \quad \quad
\quad+1883667666025x^5+383697949650x^6+41553867355x^7+2657363995x^8\\[5pt]
&\quad \quad \quad
\quad+104687375x^9+2553850x^{10}+37375x^{11}+299x^{12}+x^{13}\\[5pt]
&M_{28}(x)=2674440x+571944706335x^2+33696177453720x^3+100188554780355x^4\\[5pt]
&\quad \quad \quad
\quad+61882893062850x^5+15038660453130x^6+1925587971450x^7\\[5pt]
&\quad \quad \quad
\quad+146942256825x^8+7060951170x^9+218017800x^{10}+4296474x^{11}\\[5pt]
&\quad \quad \quad
\quad+51975x^{12}+350x^{13}+x^{14}\\[5pt]
&M_{30}(x)=9694845x+6721306583805x^2+672654700490610x^3+2722638343622385x^4\\[5pt]
&\quad \quad \quad
\quad+2089360244433195x^5+600751303879170x^6+89678011487445x^7\\[5pt]
&\quad \quad \quad
\quad+8005505867775x^8+456368933475x^9+17125516044x^{10}+426120345x^{11}\\[5pt]
&\quad \quad \quad \quad+6929405x^{12}+70470x^{13}+405x^{14}+x^{15}.
\end{align*}
It is easily verified that  $\{M_{2n,k}\}_{1\leq k\leq n}$ is
$4$-log-concave for $1\leq n \leq 15$ and $\{M_n(x)\}_{1\leq n \leq
15}$ is strongly $q$-log-convex.

If the coefficients of the $\mathcal{H}\{A_i(q)\}$ are not only
 positive but also log-concave, then we say that $\{A_i(q)\}$ is
log-concave of order two with respect to the $2$-$q$-log-convexity.
In general, we say that $\{a_i\}$ is log-concave of order $k$ with
respect to the $k$-$q$-log-convexity, if the all polynomials in
$\mathcal{H}^k\{A_i(q)\}$ are log-concave, and that $\{A_i(q)\}$ is
$q$-$\infty$-log-concave if all the polynomials in
$\mathcal{H}^k\{A_i(q)\}$ are log-concave for every $k\geq 1$.

When $1\leq n\leq 14$, we find the polynomial sequence
$\{M_{2n}(x)\}$ is $q$-$3$-log-concave. This leads us to surmise the
polynomial sequence $\{M_{2n}(x)\}$ is $q$-$\infty$-log-convex.

Finally, we recall the definition of the Boros-Moll polynomials
which are a class of the Jacobi polynomials, and are also denoted by
$P_n(a)$ as in \cite{kau}:
\begin{equation}
P_n(a)=\sum_{i=0}^n d_i(n) a^i,\end{equation}
where
\begin{align}
d_i(n)=2^{-2n}\sum_{k=i}^n2^k{2n-2k \choose n-k}{n+k \choose
k}{k\choose i}.
\end{align}

\end{document}